\documentclass[12pt]{amsart}   

\usepackage{amsmath}    
\usepackage{amsfonts}  
\usepackage{amssymb}
\usepackage{graphicx}   

\newtheorem{theorem}{Theorem}[section]
\newtheorem{corollary}[theorem]{Corollary}
\newtheorem{lemma}[theorem]{Lemma}

\begin{document}

\title{Bounding multiplicative energy by the sumset}
\author{J\'ozsef Solymosi}\thanks{Research was supported by OTKA and NSERC grants and by a Sloan Fellowship.\\ \tt solymosi@math.ubc.ca}
\date{\today}

\begin{abstract}
We prove that the sumset or the productset of any finite set of real numbers, $A,$
is at least $|A|^{4/3-\varepsilon},$ improving earlier bounds. Our main tool is a
new upper bound on the multiplicative energy, $E(A,A).$
\end{abstract}

\maketitle

\section{Introduction}
The \emph{sumset} of a finite set of an additive group, $A,$ is defined by
$$A+A=\{ a + b : a,b\in A\}.$$ 
The \emph{productset} and \emph{ratioset} are defined in a similar way.
$$AA=\{ ab : a,b\in A\},$$ and
$$A/A=\{ a/b : a,b\in A\}.$$  

A famous conjecture of Erd\H os and Szemer\'edi \cite{ESz} asserts that for any finite set of integers, $M,$
\[\max\{|M+M|,|MM|\}\geq |M|^{2-\varepsilon},\]
where $\varepsilon \rightarrow 0$ when $|M|\rightarrow \infty .$ They proved that
\[\max\{|M+M|,|MM|\}\geq |M|^{1+\delta},\]
for some $\delta >0.$ In a series of papers, lower bounds on $\delta$ were find.
$\delta \geq 1/31$ \cite{N},  $\delta \geq 1/15$ \cite{F}, $\delta \geq 1/4$ \cite{E1}, and $\delta \geq 3/14$ \cite{S2}.
The last two bonds were proved for finite sets of real numbers. 

\section{Results}
Our main result is the following.

\begin{theorem}\label{main}
Let $A$ be a finite set of positive real numbers. Then \[|AA||A+A|^2 \geq {|A|^4\over {4\lceil\log|A|\rceil}}\] holds.
\end{theorem}
The inequality is sharp - up to the power of the $\log$ term in the denominator - when $A$ is the set of the first $n$ natural numbers.
Theorem \ref{main} implies an improved bound on the sum-product problem.

\begin{corollary}
Let $A$ be a finite set of positive real numbers. Then 
\[\max\{|A+A|,|AA|\}\geq {|A|^{4/3}\over2\lceil\log|A|\rceil^{1/3}}\] holds
\end{corollary}

\subsection{Proof of Theorem \ref{main}}
To illustrate how the 
proof goes, we are making two unjustified and usually false assumptions, which
are simplifying the proof.  Readers, not interested about this "handwaving",
will find the rigorous argument about 20 lines below.  

Suppose that $AA$ and $A/A$ have
the same size,  $|AA|\approx |A/A|,$
and any element of  $A/A$ has about the same number of representations as any other. This means
that for any reals $s, t \in A/A$ the two numbers $s$ and $t$ have the same multiplicity, 
$|\{(a,b) | a, b \in A, a/b=s\}|\approx |\{(b,c) | b, c \in A, b/c=t\}|.$ A geometric interpretation
of the cardinality of $A/A$ is that the Cartesian product $A\times A$ is covered by
$|A/A|$ concurrent lines going through the origin. Label the rays from the origin
covering the points of the Cartesian product anticlockwise by 
$r_1, r_2,\ldots ,r_m,$ where $m=|A/A|.$

Our assumptions imply that each ray is incident to $|A|^2/|AA|$
points of $A\times A.$  Consider the elements of $A\times A$ as two
dimensional vectors. The sumset $(A\times A) + (A\times A) $ is the same set as $ (A+A)\times (A+A).$
We take a subset, $S,$ of this sumset, 
\[S=\bigcup_{i=1}^{m-1}(r_i \cap A\times A) + (r_{i+1}\cap A\times A)\subset (A+A)\times (A+A). \]
Simple elementary geometry shows (see the picture below) that the sumsets in the terms are disjoint and each term has \\
$|r_i \cap A\times A||r_{i+1}\cap A\times A|$ elements. Therefore
\[|S|=|AA|(|A|^2/|AA|)^2\leq |A+A|^2.\]
After rearranging the inequality we get $|A|^4\leq |AA||A+A|^2,$ as we wanted. Now we will show 
a rigorous proof based on this observation.

\medskip

We are going to use the notation of \emph{multiplicative energy}.  The name of this quantity comes from
a paper of Tao \cite{T}, however its discrete version was used earlier, like in \cite{E2}.

\medskip

Let $A$ be a finite set of reals. The multiplicative energy of $A,$ denoted by $E(A),$ is given by
\[E(A)=|\{(a,b,c,d)\in A^4 | \quad\exists \lambda \in {\bf{R}} : (a,b)=(\lambda c, \lambda d)\}|.\]

\medskip

In the notatation of Gowers \cite{G}, the quantity $E(A)$ counts the number 
of quadruples in $\log{A}.$ 

\medskip

To complete the proof of Theorem \ref{main} we show the following lemma.

\begin{lemma}\label{energy}
Let $A$ be a finite set of positive real numbers. Then \[{E(A)\over {\lceil\log{|A|}\rceil}}\leq 4|A+A|^2.\]
\end{lemma}

Theorem \ref{main} follows from Lemma \ref{energy} via the Cauchy-Schwartz type inequality
\footnote{This simple inequality appears in the literature in various places like in the proof of 
Theorem 2.4  \cite{E2}, Remark 4.2 in  \cite{T},  or  Corollary 2.10 in \cite{TV}.}   

\[E(A)\geq {|A|^4\over |AA|}.\]

\qed

\subsection{Proof of Lemma \ref{energy}}
Another way of counting $E(A)$ is the following:
\begin{equation}
E(A)=\sum_{x\in A/A} |xA\cap A|^2.
\end{equation}

The summands on the right hand side can be partitioned into $\lceil\log|A|\rceil$ classes according to the
size of $xA\cap A.$

\[E(A)=\sum_{i=0}^{\lceil\log{|A|}\rceil}\sum_{x\atop 2^i\leq|xA\cap A|<2^{i+1} } |xA\cap A|^2 \] 

There is an index, $I,$ that 

\[{E(A)\over {\lceil\log{|A|}\rceil}}\leq \sum_{x\atop 2^I\leq|xA\cap A|<2^{I+1} } |xA\cap A|^2 \] 

Let $D=\{s: 2^I\leq|sA\cap A|<2^{I+1} \},$ and let $s_1<s_2<\ldots<s_m$ denote the elements
of $D,$ labeled in increasing order. 

\begin{equation}\label{energy2}
{E(A)\over {\lceil\log{|A|}\rceil}}\leq \sum_{x\atop 2^I\leq|xA\cap A|<2^{I+1} } |xA\cap A|^2 < m2^{2I+2}.
\end{equation}

Each line $l_j : y=s_jx,$ where $1\leq j\leq m,$ is incident
to at least $2^I$ and less than $2^{I+1}$ points of $A\times A.$  For easier counting we add an extra line to the set, $l_{m+1},$ 
the vertical line through the smallest element of $A,$ denoted by $a_1.$  Line $l_{m+1}$ has $|A|$ points from $A\times A,$  however
we are considering only the orthogonal projections of the points of $l_m.$ (fig. 1) 

The sumset\footnote{As usual, 
by the sum of two points on ${\bf{R}}^2$ we mean the point which is the sum of their position vectors.}, 
$(l_i\cap A\times A) + (l_k\cap A\times A),$ $1\leq j < k\leq m,$ has size $|l_i\cap A\times A||l_k\cap A\times A|,$ 
which is between $2^{2I}$ and $2^{2I+2}.$ Also, the sumsets along consecutive line pairs are
disjoint, i.e. 
\[((l_i\cap A\times A) + (l_{i+1}\cap A\times A))\cap((l_k\cap A\times A) + (l_{k+1}\cap A\times A))=\emptyset ,\]
for any  $1\leq j < k \leq m.$

\begin{figure}[htbp]
\begin{center}
\includegraphics[scale=0.37]{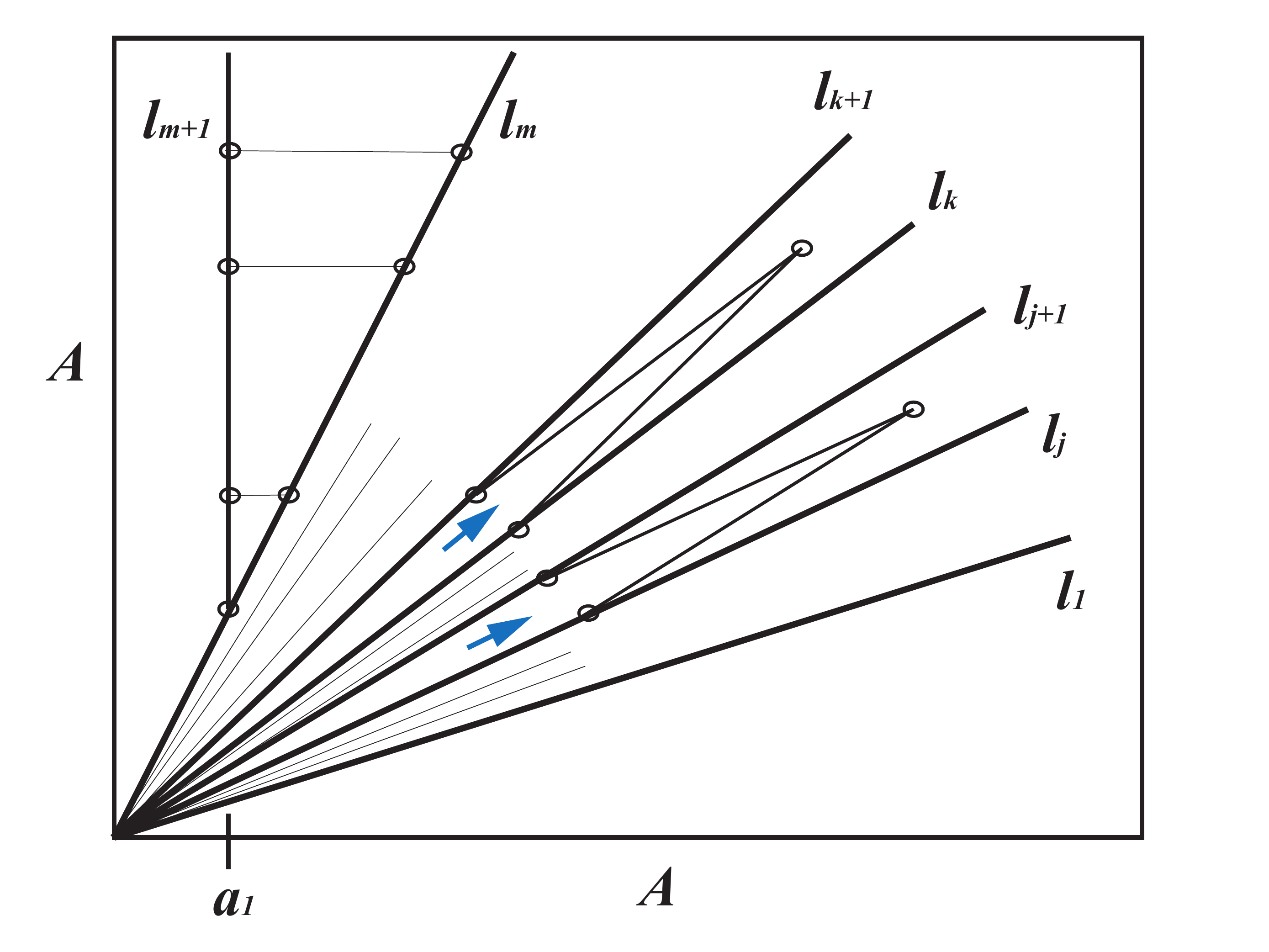}
\caption{}
\label{default}
\end{center}
\end{figure}

The sums are elements of $(A+A) \times (A+A),$ so we have the following inequality.

\[m2^{2I}\leq\left|\bigcup_{i=1}^{m}(l_i \cap A\times A) + (l_{i+1}\cap A\times A)\right|\leq |A+A|^2. \]
The inequality above with inequality (\ref{energy2}) proves the lemma.

\qed

\subsection{Remarks}
Let $A$ and $B$ be finite sets of reals. The multiplicative energy, $E(A,B),$ is given by
\[E(A,B)=|\{(a,b,c,d)\in A\times B\times A\times B | \quad\exists \lambda \in {\bf{R}} : (a,b)=(\lambda c, \lambda d)\}|.\]
In the proof of Lemma \ref{energy} we did not use the fact that $A=B,$ the proof works for the
asymmetric case as well. Suppose that $|A|\geq |B|.$ With the lower bound on the multiplicative energy

\[E(A,B)\geq {|A|^2|B|^2\over |AB|}\]

our proof gives the more general inequality
\[{|A|^2|B|^2\over {|AB|}}\leq 4 {\left\lceil\log{|B|}\right\rceil}|A+A||B+B|.\]

\section{Very small productsets}
In this section we extend our method from two to higher dimensions. We are going to consider
lines though the origin as before, however there is no notion of consecutiveness among these
lines in higher dimensions available. We will consider them as points in the projective real space and
will find a triangulation of the pointset. The simplices of the triangulation will define the
neighbors among the selected lines.
\medskip

The sum-product bound in Theorem \ref{main} is asymmetric. It shows that the productset should be very large
if the sumset is small. On the other hand it says almost nothing in the range where the productset is
small. For integers, Chang \cite{MCC2} proved that there is a function $\delta(\varepsilon)$ that if $|AA|\leq |A|^{1+\varepsilon}$
then $|A+A|\geq |A|^{2-\delta},$ where $\delta \rightarrow 0$ if $\varepsilon \rightarrow 0.$ Similar result
is not known for reals. It follows from Elekes' bound \cite{E1} (and also from Theorem \ref{main}) that 
there is a function $\delta(\varepsilon)$ that if $|AA|\leq |A|^{1+\varepsilon}$
then $|A+A|\geq |A|^{3/2-\delta},$ where $\delta \rightarrow 0$ if $\varepsilon \rightarrow 0.$ 
We prove here a generalization of this bound for $k$-fold sumsets. For any integer $k\geq 2$ the $k$-fold subset of $A,$
denoted by $kA$ is the set
\[kA=\{a_1+a_2+\ldots +a_k  | a_1,\ldots ,a_k \in A\}.\]

\begin{theorem}\label{kfold}
For any integer $k\geq 2$ there is a function $\delta=\delta_k(\varepsilon)$ that if $|AA|\leq |A|^{1+\varepsilon}$
then $|kA|\geq |A|^{2-1/k-\delta},$ where $\delta \rightarrow 0$ if $\varepsilon \rightarrow 0.$ 
\end{theorem}

\begin{proof}
We can suppose that  $A$ has only positive elements WLOG. Let $|AA|\leq |A|^{1+\varepsilon}.$ By a Pl\"unnecke type inequality 
(Corollary 5.2 \cite{Ru} or Chapter 6.5 \cite{TV}) we
have $|A/A|\leq |A|^{1+2\varepsilon}.$ Consider the $k$-fold Cartesian product $A\times A\times \ldots \times A,$
denoted by $\times^kA.$ It can be covered by no more than $|A/A|^{k-1}$ lines going through the origin.
\begin{figure}[htbp]
\begin{center}
\includegraphics[scale=0.37]{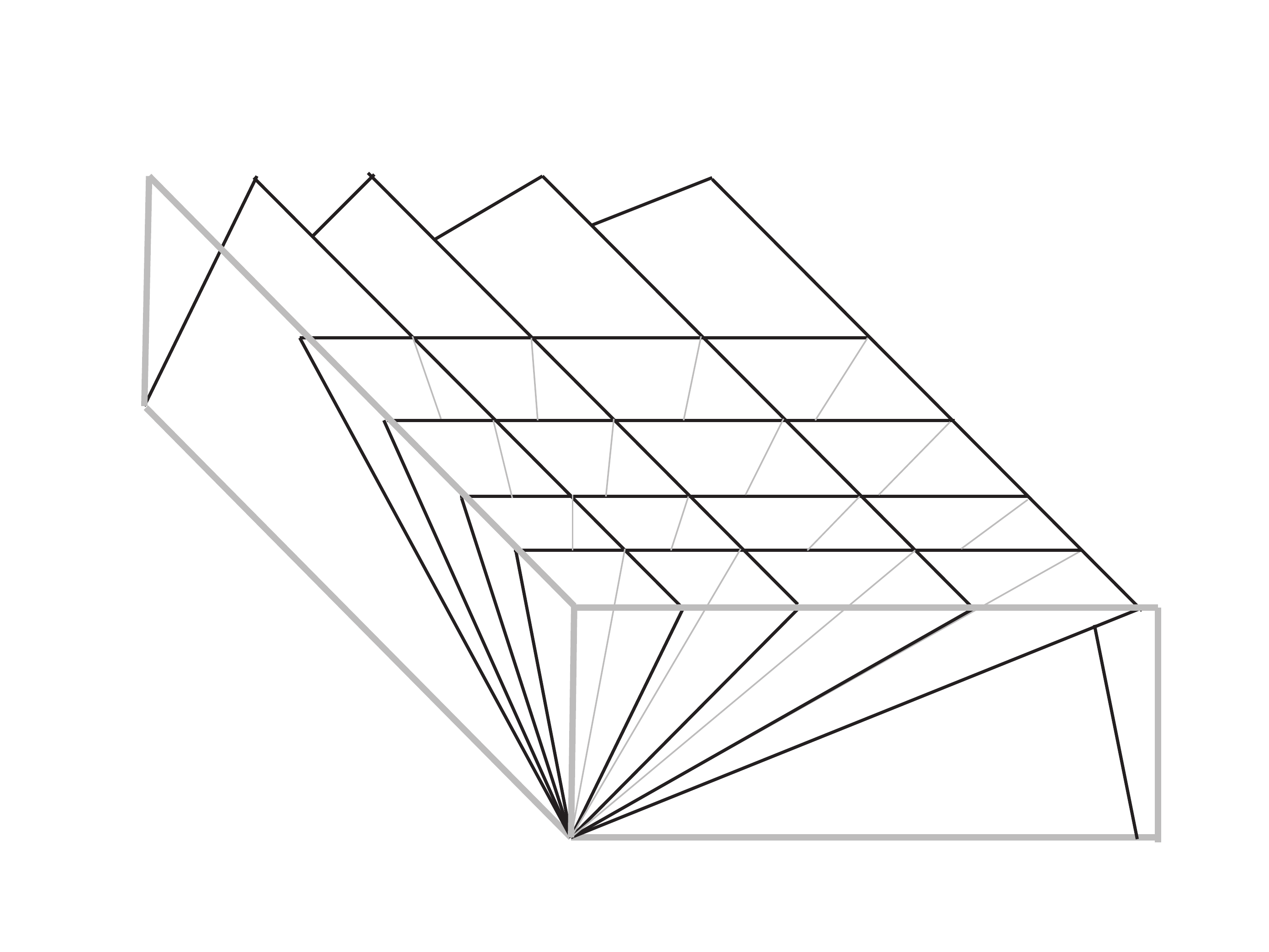}
\caption{}
\label{grid}
\end{center}
\end{figure}
Picture \ref{grid} illustrates the $k=3$ case. Let $H$ denotes the set of lines through the origin containing at least
$|A|^{1-2\varepsilon(k-1)}/2$ points of $\times^kA.$ 
With this selection, the lines in $H$ cover at least half of the points in $\times^kA$ since 
\[{|A|^{1-2\varepsilon(k-1)}\over 2}|A/A|^{k-1}={|A|^k\over {2|A|^{(1+2\varepsilon)(k-1)}}}|A/A|^{k-1}\leq {|A|^k\over 2}.\]
As no line has more than $|A|$ points common with $\times^kA,$ therefore $|H|\geq |A|^{k-1}/2.$
The set of lines, $H,$ represents a set of points, $P,$  in the projective real space ${\bf RP}^{k-1}.$ Point set $P$ has full dimension
$k-1$ as it has a nice symmetry. The symmetry follows from the Cartesian product structure; if a point with 
coordinates $(a_1, \ldots, a_k)$ is in $P$ then the point $(\sigma(a_1),\ldots, \sigma(a_k))$ is also in 
$P$ for any permutation $\sigma\in S_k.$ Let us triangulate $P.$  By triangulation we mean a 
decomposition of the convex hull of $P$ into non-degenarate, $k-1$ 
dimensional, simplices such that the intersection of any two is the empty set or 
a face of both simplices and the vertex set of the triangulation  is  $P$. It is not obvious that such triangulation 
always exists. For the proof we refer to Chapter 7 in \cite{GKZ} or
Chapter 2 in \cite{LRS}.    Let ${\tau}(P)$ be a triangulation of $P.$ We say that $k$ lines 
$l_1,\ldots,l_k\in H$ form a simplex if the corresponding points in $P$ are vertices of a simplex
of the triangulation. We use the following notation for this: $\{l_1,\ldots,l_k\} \in \tau(P).$
In the two-dimensional case we used that the sumsets of points on consecutive lines
are disjoint. Here we are using that the interiors of the simplices are disjoint, therefore
sumsets of lines of simplices are also disjoint. Note that we assumed that $A$ is positive,
so we are considering convex combinations of vectors with positive coefficients.
Let  $\{l_1,\ldots,l_k\} \in \tau(P)$ and  $\{l'_1,\ldots,l'_k\} \in \tau(P)$ are two distinct
simplices. Then
\[\left(\sum_{i=1}^kl_i\cap\times^kA\right)\bigcap \left(\sum_{i=1}^kl'_i\cap\times^kA\right) =\emptyset .\] 
Also, since the $k$ vectors parallel to the lines $\{l_1,\ldots,l_k\} \in \tau(P)$ are linearly independent, all sums are distinct,
\[\left|\sum_{i=1}^kl_i\cap\times^kA\right|=\prod_{i=1}^k\left|l_i\cap\times^kA\right|.\]
Now we are ready to put everything together into a sequence of inequalities proving Theorem \ref{kfold}.
\[|kA|^k\geq \sum_{\{l_1,\ldots,l_k\} \in \tau(P)}\left|\sum_{i=1}^kl_i\cap\times^kA\right|\geq {|A|^{k-1}\over 2k}\prod_{i=1}^k\left|l_i\cap\times^kA\right|.\]
Every line is is incident to at least $|A|^{1-2\varepsilon(k-1)}/2$ points of $\times^kA,$ therefore
\[|kA|^k\geq  {|A|^{k-1+k(1-2\varepsilon(k-1))}\over 2k2^k}={|A|^{2k-1-2k(k-1)\varepsilon}\over k2^{k+1}}.\]
Taking the $k$-th root of both sides we get the result we wanted to show
\[|kA|\geq c_k|A|^{2-1/k - 2(k-1)\varepsilon}.\]
\end{proof}

\end{document}